\documentclass[12pt,leqno,fleqn]{amsart} 

\usepackage{amsmath,amstext,amsthm,amssymb,amsxtra} 
\usepackage[color=green!15,bordercolor=red,linecolor=red!30]{todonotes}  
\usepackage{ulem}

\usepackage{hyperref} 
\hypersetup{
    colorlinks=true,       
    linkcolor=blue,          
    citecolor=magenta,        
    filecolor=magenta,      
    urlcolor=cyan,          
}

\newtheorem{theorem}{Theorem}[section]
\newtheorem{lemma}[theorem]{Lemma}

\newtheorem{proposition}[theorem]{Proposition}

\theoremstyle{remark}

\def\C{{\mbox{\rm\kern.24em
\vrule width.03em height1.43ex depth-.052ex \kern-.26em C}}}
\def\QSet{\mbox{\rm\kern.24em
\vrule width.03em height1.48ex depth-.051ex \kern-.26em Q}}
\def\R{{\mbox{\rm I\kern-.22em R}}}
\def\E{{\mbox{\rm I\kern-.22em E}}}
\def\P{{\bf P}}

\def\K{{\bf K}}

\def\Re{{\rm Re}}

\def\J{{\bf J}}

\def\\P{{\mathcal P}}

\def\bas{\begin{align*}}
\def\eas{\end{align*}}
\def\bi{\begin{itemize}}
\def\ei{\end{itemize}}
\def\emph#1{{\it #1}}

\begin{document}
\title{The Hilbert transform does not map $L^1(Mw)$ to $L^{1,\infty}(w)$}

\author{Maria Carmen Reguera}
\address{School of Mathematics, Georgia Institute of Technology, Atlanta GA 
30332}
\email{mreguera@@math.gatech.edu}

\author{Christoph Thiele}
\address{Department of Mathematics, UCLA, Los Angeles CA 90095-1555}
\email{thiele@@math.ucla.edu}


\thanks{M.C.R. partially supported by NSF grant DMS 0968499.} 
\thanks{C.Th. partially supported by NSF grant DMS 1001535.}
\date{\today}

\begin{abstract}

We disprove the following a priori estimate
for the Hilbert transform $H$ and the Hardy Littlewood maximal
operator $M$:
$$\sup_{t>0}t w\{x\in \R:|Hf(x)|>t\}\le C\int |f(x)| Mw(x) \, dx\ \ .$$
This is a sequel to paper \cite{reguera} by the first author, 
which shows the existence of a Haar multiplier operator for which the inequality holds.

\end{abstract}

\maketitle

\section{Introduction and statement of main result}

In \cite{FS}, C. Fefferman and E. Stein observed the following
a priori estimate for the Hardy Littlewood maximal operator $M$:
$$\sup_{t>0}t\,  w\{x\in \R:|Mf(x)|>t\}\le C\int |f(x)| Mw(x) \, dx\ \ .$$
Here the weight $w$ is a non-negative, locally integrable function, and $w(E)$ 
denotes the integral of the weight over the set $E$. We give a negative answer to 
the question whether such an inequality holds
when the Hardy Littlewood maximal operator on the left hand side is replaced 
by the Hilbert transform. For a discussion of the history of this question 
we refer to \cite{reguera}.

\begin{theorem}\label{main}
For each constant $C>0$ there is a weight function
$w$ on the real line and an integrable compactly supported function $f$ 
and a $t>0$ such that 
$$ t\, w\{x\in \R:|Hf(x)|>t\}\ge C\int |f(x)| Mw(x) \, dx\ \ .$$
\end{theorem}
Similarly as in \cite{reguera}, we prove 
Theorem \ref{main} as a consequence of the following:
\begin{proposition}\label{cuperez}
For each constant $C>0$ there is an everywhere positive weight function
$w$ on the real line and an integrable compactly supported function $f$ 
and a $t>0$ such that 
\begin{equation}\label{weak2}
t^2\,  {w}\{x\in\R: |Hf(x)|>t\}\ge C\int  |f(x)|^2 \left(\frac{Mw(x)}{w(x)}\right)^2 w(x)\, dx  \ \ .
\end{equation}
\end{proposition}
The reduction to Proposition \ref{cuperez} is taken
from \cite{cuperez}, we sketch the argument at the end of this paper.
Following \cite{reguera} further, we reduce 
Proposition \ref{cuperez} to the dual proposition:
\begin{proposition}\label{dualcp}
For each constant $C$ there is a nontrivial weight $w$ on the real line such that
$$\|H(w1_{[0,1)})\|_{L^2(w/(Mw)^2)}
\ge C \|1_{[0,1)}\|_{L^2(w)}
\ \ .$$
\end{proposition}
Our construction of the weight $w$ is a somewhat simpler variant of 
the construction in \cite{reguera}. It was discovered during
a stimulating summer school on ``Weighted estimates for
singular integrals'' at Lake Arrowhead, Oct 3-8. 2010.

\section{Proof of Theorem \ref{cuperez}}

Recall that a triadic interval $I$ is of the form $[3^{j}n, 3^j(n+1))$ with integers $j,n$.
Denote by ${I^m}$ the triadic interval of one third the length of $I$ which contains the center of $I$.

Fix an integer $k$ which will be chosen large enough depending on the constant $C$ in Proposition 
\ref{dualcp}. Define $\K_0$ to be $\{[0,1)\}$ and recursively for $i\ge 1$:
$$\J_i:=\{K^m: K\in \K_{i-1}\}\ \ ,$$
$$\K_i:=\{K : K {\rm \ triadic},\  |K|=3^{-ik},\  K\subset \bigcup_{J\in\J_i}J\}\ \ .$$

Proceeding recursively from the larger to the smaller intervals, 
we choose for each $J\in \J:=\bigcup_{i\ge 1} \J_i$ a sign $\epsilon(J)\in \{-1,1\}$.
More precisely, $\epsilon(J)$ depends on the values $\epsilon(J')$ with $|J'|>|J|$. 
The exact choice will be specified below.
Define for each $J\in\J$ the interval
$I(J)$ to be the triadic interval of length $3^{1-k}|J|$ whose right endpoint
equals the left endpoint of $J$ if $\epsilon(J)=1$, and whose
left endpoint equals the right endpoint of $J$ if $\epsilon(J)=-1$.
Note that $I(J)$ has the same length as the intervals in $\K_i$.

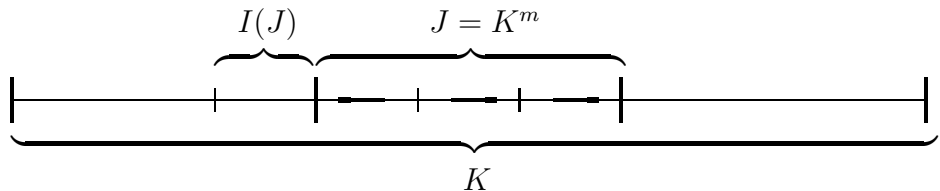
\setlength{\unitlength}{1.5mm}
\begin{picture}(100,23)



\put(37,13){$\overbrace{\rule{4.1cm}{0cm}}$}
\put(10,7){$\underbrace{\rule{12.3cm}{0cm}}$}
\put(28,13){$\overbrace{\rule{1.3cm}{0cm}}$}

\put(50,2){$K$}
\put(47,16){$J=K^m$}
\put(30,16){$I(J)$}

\put(10,10){\line(1,0){81}}
\multiput(28,9)(9,0){5}{\line(0,1){2}}
\linethickness{1pt}
\multiput(10,8)(27,0){4}{\line(0,1){4}}
\multiput(40,10)(9,0){3}{\line(1,0){3}}
\linethickness{2pt}
\put(39,10){\line(1,0){1}}
\put(52,10){\line(1,0){1}}
\put(61,10){\line(1,0){1}}

\end{picture}

Next we define a sequence of absolutely continuous measures on $[0,1]$.
We continue to use the same symbol for a measure and its Lebesgue density.
Let $w_0$ be the uniform measure on $[0,1)^m\cup I([0,1))^m$
with total mass $1$. Recursively we define the measure $w_i$ by
the following properties: It coincides with $w_{i-1}$ on the complement of 
$\bigcup_{K\in \K_{i}}K$.
For $K\in \K_i$ we have $w_i(K)=w_{i-1}(K)$
and the restriction of $w_i$ to $K$ is supported and uniformly distributed on 
$K^m\cup I(K^m)$.

Let $w$ be the weak limit of the sequence $w_i$ and note 
that $w$ is supported on $\bigcup_{J\in \J} I(J)$.
For $K\in \K_i$, $J\in \J_{i}$, $x\in I(J)$, and any triadic
interval $K'$ with $|K'|\ge |K|$ we have
\begin{equation}\label{intcompare}
w(x)=\frac{w(I(J))}{|I(J)|}=\frac{w(K)}{|K|}\ge \frac{w(K')}{|K'|}\ \ .
\end{equation}
We claim that for $J\in \J$ and $x\in I(J)^m$ we have
\begin{equation}\label{mwcompare}
Mw(x)\le 7w(x)\ \ .
\end{equation}
To see this, let $I$ be a (not necessarily triadic) interval containing $x$. 
If $I$ is contained in $I(J)$, then by the first identity of (\ref{intcompare}) the 
average of $w$ over $I$ equals $w(x)$. If $I$
is not contained in $I(J)$, then $|I|\ge |I(J)|/3$.
Let $\K'$ be the collection of triadic intervals of
length $|I(J)|$ which intersect $I$ and note that
$$\sum_{K'\in \K'} |K'|\le |I|+2|I(J)|\le 7|I|$$
because at most two intervals in $\K'$ are not entirely covered by $I$.
With (\ref{intcompare}) we conclude that the average of $w$ over $I$ is
no more than $7w(x)$, which completes the proof of (\ref{mwcompare}).

\begin{lemma}\label{hestimate}
For $K\in \K_{i}$, $J=K^m$, $x\in I(J)^m$, and $k>3000$ we have
$$|Hw(x)|\ge (k/3) w(x)\ \ .$$
\end{lemma}
This Lemma 
proves Proposition \ref{dualcp},
because with (\ref{mwcompare}) and since 
$w$ is constant on every $I(J)$
we have
$$49 \|H w\|_{L^2(w/(Mw)^2)}^2\ge (k^2/9) \sum_{J\in \J}
\int_{I(J)^m}w(y)\, dy \ge (k^2/27) \|1_{[0,1)}\|_{L^2(w)}^2 \ \ .$$
Proof of Lemma \ref{hestimate}: We split the principal value
integral for $Hw(x)$ into six summands:
\begin{equation}\label{a1}
p.v. \int_{I(J)} \frac {w(y)}{y-x} \, dy 
\end{equation}
\begin{equation}\label{a2}
 + \int_{J} \frac {w(y)}{y-x} \, dy 
\end{equation}
\begin{equation}\label{a3}
 + \int_{K^c} (\frac {w(y)}{y-x}-\frac {w(y)}{y-c(J)})\, dy
\end{equation}

\begin{equation}\label{a4}
 +  \int_{(\bigcup_{\K_i} K')^c} \frac {w(y)}{y-c(J)} \, dy
\end{equation}

\begin{equation}\label{a5}
+\sum_{K'\in \K_i\setminus \{K\}} \int_{K'}\frac{w(y)}{y-c(J)}- \frac {w(y)}{c(K')-c(J)} \, dy
\end{equation}

\begin{equation}\label{a6}
+\sum_{K'\in \K_i\setminus\{K\}} \int_{K'} \frac {w(y)}{c(K')-c(J)} \, dy\ \ .
\end{equation}

The terms (\ref{a4}) and (\ref{a6}) remain unchanged if we replace
$w$ by $w_{i}$ and hence depend only on the
choices of $\epsilon(J')$ with $|J'|>|J|$. The integrand of $(\ref{a2})$ is positive 
or negative depending on $\epsilon(J)$.
Specify the choice of $\epsilon(J)$ so that the sign of (\ref{a2}) 
equals the sign of (\ref{a4})+(\ref{a6}). If the latter is zero,
we may arbitrarily set $\epsilon(J)=1$. We estimate 
$$|(\ref{a2})|\ge 
\sum_{K'\in \K_{i+1}, K'\subset J} \int_{K'} \frac {w(y)}{|y-x|} \, dy  $$
$$\ge 
\sum_{K'\in \K_{i+1}, K'\subset J} \frac {w(K')}{\sup_{y\in K'} |y-x|}  $$
$$\ge 
\sum_{n=1}^{3^k} \frac {1}{n+1} \frac{w(I(J))}{|I(J)|} 
\ge {(k/2)} w(x)\ \ . $$
The remaining terms are small error terms, we estimate with $\delta=|I(J)^m|$
:
$$
|(\ref{a1})|= | \int_{I(J)\setminus [x-\delta,x+\delta]} \frac {w(y)}{y-x} \, dy |
\le 3 w(x) \ \ ,$$
$$
 |(\ref{a3})|
\le 4  \sum_{|K'|=|K|, K'\neq K}
\int_{K'} \frac {|x-c(J)|}{|y-c(J)|^2} w(y)\, dy
$$
$$
\le 8  \sum_{|K'|=|K|, K'\neq K} \frac {|x-c(J)|}{|c(K')-c(J)|^2} w(K')
$$
$$
\le 16  \sum_{n=1}^\infty \frac {1}{(n-3/4)^2}
\frac{w(I(J))}{|I(J)|} \le 200 w(x)\ \ ,$$
$$
|(\ref{a5})|\le 
 4 \sum_{K'\in \K_i} \int_{K'}\frac{|y-c(K')|}{|c(K')-c(J)|^2} w(y)\, dy\ \ ,
$$
and the last expression is dominated by the same final bound as (\ref{a3}).
Putting all estimates together, we have 
$$|(\ref{a1})+(\ref{a2})+(\ref{a3})+(\ref{a4})+(\ref{a5})+(\ref{a6})|$$
$$\ge |(\ref{a2})+(\ref{a4})+(\ref{a6})| - |(\ref{a1})|- |(\ref{a3})|- |(\ref{a5})|$$
$$\ge |(\ref{a2})| - |(\ref{a1})|- |(\ref{a3})|- |(\ref{a5})|$$
$$\ge (k/2-403) w(x)\ \ .$$
This completes the proof of Lemma \ref{hestimate} and thus Theorem \ref{dualcp}.

\section{Remarks}

\subsection{More general kernels} The construction can be generalized to 
apply to more general kernels, including those with even symmetry, such as for example $\Re(|x|^{-1+\alpha i})$ with $\alpha\neq 0$. Choose $J$ to be the union of $3^{k-1}$ not necessarily 
adjacent but appropriately chosen intervals of length $3^{-k}|K|$ contained in $K$,  and $I(J)$ an appropriate 
further interval of this length well inside $K$, so that the kernel of the 
Calderon Zygmund operator for $x\in I(J)^m$ has sufficient positive or 
negative bias on $J$.

\subsection{Weights in Theorem \ref{main}} 
We specify weights satisfying Theorem (\ref{main}).
Fix a constant $C$ as in Proposition (\ref{dualcp}) and consider $k$ and the weight $w$ 
constructed above. We slightly change $w$ to make it positive by adding $ce^{-x^2}$
for sufficiently small $c$ so as to not change the conclusion of Proposition 
(\ref{dualcp}). We may normalize the measure to be probability measure and call
the remaining measure $w$ again. The conclusion of Proposition \ref{dualcp} can be written:
\begin{equation}\label{explicitdualpc}
(\int (Hw(x))^2 \frac{w(x)}{(Mw(x))^2}\, dx)^{1/2} \ge C\ \ .
\end{equation}
Multiplying both sides of (\ref{explicitdualpc})
by the left hand side of (\ref{explicitdualpc}),
setting $f= (Hw)w/(Mw)^2$ and using essential self-duality of $H$ we obtain
\begin{equation}\label{dualpctopc}
 | \int w(x) Hf(x) \, dx|\ge C(\int f(x)^2 \frac{(Mw(x))^2}{(w(x))^2} w(x)\, dx)^{1/2}\ \ .
\end{equation}
Letting $f^*$ be the non-increasing rearrangement of $Hf$ on $[0,1]$, we may estimate the 
left hand side of (\ref{dualpctopc})
$$\int_0^1 f^*(y)\, dy\le 2\sup_{y\in[0,1]} y^{1/2}f^*(y)=2\sup_{t>0}  w(\{x:|Hf(x)|\ge t\})^{1/2} t\ \ .$$
Hence Proposition \ref{cuperez} holds for the constant $C/2$ with the weight $w$
and some existentially chosen $t$. Now let $E$ be the set on the left hand side of Proposition 
\ref{cuperez} for the given $w$, $f$, and appropriate $t$, then we have
$$M(w1_E)(x)=\sup_{x\in I}\frac{\int_I w}{\int_I 1}\frac{\int_I 1_E w}{\int_I w}\le Mw(x) M_w 1_E(x)\ \ ,$$
where $M_w$ denotes the Hardy Littlewood maximal function with respect to the weight $w$.
With H\"older's inequality we obtain
$$\int |f(x)| M(w1_E )(x) \, dx 
\le \left(\int |f(x)|^2 \frac{Mw(x)^2} {w(x)}\, dx\right)^{1/2}  \|M_w 1_E\|_{L^2(w)}\ \ .$$
With the Hardy Littlewood maximal theorem with respect to the weight $w$ we can estimate
$\|M_w 1_E\|_{L^2(w)}$ by $w(E)^{1/2}$. This shows that Theorem \ref{main} holds for
the weight $w1_E$.

\subsection{A1 weights}
It remains open to date whether the a priori inequality
\begin{equation}\label{wa1}
 t\, w\{x\in \R:|Hf(x)|>t\}\le C\|w\|_{A_1}\int |f(x)| w(x) \, dx 
\end{equation}
holds, where the $A_1$ constant is defined as $\|w\|_{A_1}:=\|Mw/w\|_\infty$. 
Our construction in this paper does not seem to
address this question.  The recent preprint \cite{nrvv} has announced that the analogue of (\ref{wa1}) for
Haar multipliers is false. In  \cite{lop}, a version of (\ref{wa1}) 
has been proved with an additional logarithmic factor in the $A_1$ constant of the weight.

\end{document}